\documentclass[a4paper,12pt]{article}

\input{mymacros.sty}

\title{Dimer models for parallelograms}
\author{Kazushi Ueda and Masahito Yamazaki}
\date{}
\begin{document}

\maketitle

%
%-------------------- abstract -------------------------
%

\begin{abstract}
We discuss the relation
between dimer models and coamoebas
associated with lattice parallelograms.
We also discuss
homological mirror symmetry
for $\bP^1 \times \bP^1$,
emphasizing the role of a non-isoradial dimer model.
\end{abstract}

%
%=------------------ text starts -----------------------
%

\section{Introduction}
This is a short companion paper to
\cite{Ueda-Yamazaki_NBTMQ, Ueda-Yamazaki_toricdP},
where we discuss the relation between
dimer models, coamoebas and homological mirror symmetry,
following the work of string theorists
\cite{Feng-He-Kennaway-Vafa,
Franco-Hanany-Martelli-Sparks-Vegh-Wecht_GTTGBT,
Franco-Hanany-Vegh-Wecht-Kennaway_BDQGT,
Hanany-Herzog-Vegh_BTEC,
Hanany-Kennaway_DMTD,
Hanany-Vegh}.
We refer the readers to \cite{Ueda-Yamazaki_NBTMQ, Ueda-Yamazaki_toricdP}
and references therein
for notations and backgrounds.
We show the following in this paper:
\begin{itemize}
 \item
The linear Hanany-Vegh algorithm
produces a unique dimer model from a lattice parallelogram.
 \item
For a suitable choice of a Laurent polynomial $W$,
the zero locus $W^{-1}(0)$
behaves nicely under the argument map
$
 \Arg : (\bCx)^2 \to (\bR / \bZ)^2.
$
%$$
%\begin{array}{cccc}
% \Arg : & (\bCx)^2 & \to & T = (\bR / \bZ )^2 \\
% & \vin & & \vin \\
% & (x, y) & \mapsto & {\displaystyle \frac{1}{2 \pi}(\arg(x), \arg(y))}.
%\end{array}
%$$
 \item
%We prove \cite[Conjecture 6.2]{Ueda-Yamazaki_toricdP}
%for $\bP^1 \times \bP^1$,
%which states that
The dimer model is the argument projection
of a graph on $W^{-1}(0)$,
which encodes Floer-theoretic information
of vanishing cycles.
\end{itemize}
As a corollary,
one obtains a torus-equivariant generalization
of homological mirror symmetry for $\bP^1 \times \bP^1$
proved by Seidel \cite{Seidel_VC2}.
We also discuss the role of a non-isoradial dimer model
which cannot be obtained by the linear Hanany-Vegh algorithm.

\section{Coamoebas for lattice parallelograms}
We prove the following in this section:

\begin{theorem}
For any lattice parallelogram,
there exists a Laurent polynomial $W$
satisfying the following:
\begin{itemize}
 \item
The Newton polygon of $W$ coincides with the lattice parallelogram.
 \item
The set $\scA$ of asymptotic boundaries
of the coamoeba of $W^{-1}(0)$
is the unique admissible arrangement of oriented lines
associated with the lattice parallelogram.
 \item
The coamoeba is the union of colored cells and vertices of $\scA$.
 \item
The restriction of the argument map
to the inverse image of a colored cell
is a diffeomorphism.
It is orientation-preserving if the cell is white,
and reversing if the cell is black.
 \item
The inverse image of a vertex of $\scA$
by the argument map
is homeomorphic to an open interval.
\end{itemize}
\end{theorem}

\begin{proof}
\begin{figure}[htbp]
\begin{tabular}{cc}
\begin{minipage}{.5 \linewidth}
\centering
\input{Deltaone_polygon.pst}
\caption{A lattice square}
\label{fg:Deltaone_polygon}
\end{minipage}
&
\begin{minipage}{.5 \linewidth}
\centering
\input{Deltaone_surface.pst}
\caption{The projection of $W^{-1}(0)$ to the $x$-plane}
\label{fg:Deltaone_surface}
\end{minipage}
\\
\begin{minipage}{.5 \linewidth}
\centering
\input{Deltaone_coamoeba.pst}
\caption{The coamoeba}
\label{fg:Deltaone_coamoeba}
\end{minipage}
&
\begin{minipage}{.5 \linewidth}
\centering
\input{Deltaone_lines.pst}
\caption{The dimer model}
\label{fg:Deltaone_graph}
\end{minipage}
\end{tabular}
\end{figure}
Since any lattice parallelogram
is obtained from the convex hull of
$(0, 0)$, $(0, 1)$, $(1, 0)$ and $(1, 1)$
shown in Figure \ref{fg:Deltaone_polygon}
by a translation and an integral linear transformation,
it suffices to discuss the case
$$
 W(x, y) = x y + x - y + 1.
$$
The projection
$$
\begin{array}{ccc}
 W^{-1}(0) & \to & \bCx \\
 \vin & & \vin \\
 (x, y) & \mapsto & x
\end{array}
$$
to the $x$-plane is injective,
whose image is obtained by gluing
the upper half plane
$$
 D_1 = \{ x \in \bC \mid \Im(x) > 0 \}
$$
and the lower half plane
$$
 D_2 = \{ x \in \bC \mid \Im(x) < 0 \}
$$
along four intervals
$$
 I_1 = (-\infty, -1), \ 
 I_2 = (-1, 0), \ 
 I_3 = (0, 1), \text{ and }
 I_4 = (1, \infty).
$$
as shown in Figure \ref{fg:Deltaone_surface}.
Set
\begin{equation*}
\begin{split}
 U_1 &=
 \lc (\theta, \phi) \in T \, \left| \ 
      0 < \theta < \frac{1}{2}, \ 
      0 < \phi < \frac{1}{2}
 \right. \rc, \\
 U_2 &=
 \lc (\theta, \phi) \in T \, \left| \ 
      -\frac{1}{2} < \theta < 0, \ 
      -\frac{1}{2} < \phi < 0,
 \right. \rc.
\end{split}
\end{equation*}
Then the argument map
gives an orientation-reversing homeomorphism
from $D_1$ to $U_1$ and
an orientation-preserving homeomorphism
from $D_2$ to $U_2$.
The line segments
$I_1$, $I_2$, $I_3$ and $I_4$
on the boundary of $D_1$ and $D_2$ are mapped to four points
$(\frac{1}{2}, \frac{1}{2})$,
$(\frac{1}{2}, 0)$,
$(0, 0)$ and
$(0, \frac{1}{2})$
respectively.
This shows that
the coamoeba of $W^{-1}(0)$ is
as in Figure \ref{fg:Deltaone_coamoeba}.
Its asymptotic boundaries and
the corresponding dimer model is
shown in Figure \ref{fg:Deltaone_graph}.

The uniqueness of the output of the linear Hanany-Vegh algorithm
can be shown in just the same way as the case of triangles
in \cite[Section 5.6]{Ueda-Yamazaki_NBTMQ},
and we omit the detail here.
\end{proof}

\section{Vanishing cycles for the mirror of $\bP^1 \times \bP^1$}

We prove \cite[Conjecture 6.2]{Ueda-Yamazaki_toricdP}
for $\bP^1 \times \bP^1$ in this section.
\begin{figure}[htbp]
\begin{tabular}{cc}
\begin{minipage}{.5 \linewidth}
\centering
\input{P1P1_vp1.pst}
\caption{A distinguished set of vanishing paths}
\label{fg:P1P1_vp1}
\end{minipage}
&
\begin{minipage}{.5 \linewidth}
\centering
\input{s-plane.pst}
\caption{Branch points and cuts on the $s$-plane}
\label{fg:s-plane}
\end{minipage}
\\
\ \vspace{0mm}\\
\begin{minipage}{.5 \linewidth}
\centering
\input{P1P1_surface.pst}
\caption{The glued surface $W^{-1}(0)$}
\label{fg:P1P1_surface}
\end{minipage}
&
\begin{minipage}{.5 \linewidth}
\centering
\input{s-plane_motion.pst}
\caption{Behavior of the branch points along $c_1$}
\label{fg:s-plane_motion}
\end{minipage}
% \begin{minipage}{.5 \linewidth}
% \input{s-plane_c1.pst}
% \caption{The image of $C_1$ by $\pi_0$}
% \label{fg:s-plane_c1}
% \end{minipage}
\end{tabular}
\end{figure}
\begin{figure}
\begin{minipage}{.5 \linewidth}
\centering
\input{s-plane_vc.pst}
\caption{The images of four vanishing cycles
on the $s$-plane}
\label{fg:s-plane_vc}
\end{minipage}
% \begin{minipage}{.5 \linewidth}
% \input{P1P1_surface_C1.pst}
% \caption{The vanishing cycle $C_1$}
% \label{fg:P1P1_surface_C1}
% \end{minipage}
\begin{minipage}{.5 \linewidth}
\centering
\input{P1P1_surface_cycles.pst}
\caption{Vanishing cycles on $W^{-1}(0)$}
\label{fg:P1P1_surface_cycles}
\end{minipage}
\end{figure}
The mirror of $\bP^1 \times \bP^1$ is given by the Laurent polynomial
$$
 W(x, y) = x - \frac{1}{x} + y + \frac{1}{y},
$$
whose critical points are given by
$$
 (x, y)
   = (- \sqrt{-1}, -1),
     ( \sqrt{-1}, -1),
     ( \sqrt{-1}, 1),
     (- \sqrt{-1}, 1),
$$
with critical values
$$
 (p_1, p_2, p_3, p_4)
  = ( - 2 - 2 \sqrt{-1},
      - 2 + 2 \sqrt{-1},
      2 + 2 \sqrt{-1},
      2 - 2 \sqrt{-1}).
$$
Let $(c_i)_{i=1}^4$ be a distinguished set of vanishing paths,
obtained as the straight line segments
from the origin to the critical values of $W$
as shown in Figure \ref{fg:P1P1_vp1}.
The fiber of $W$ can be realized
as a branched cover of $\bC$ by
$$
\begin{array}{cccc}
 \pi_t : & W^{-1}(t) & \to & \bC \\
 & \rotatebox{90}{$\in$} & & \rotatebox{90}{$\in$} \\
 & (x, y) & \mapsto & x - \displaystyle{\frac{1}{x}}
\end{array}
$$
for $t \in \bC$.
The fiber of $\pi_t$ at a general point $s \in \bC$
consists of four points,
which is the product of
the set of solutions of
\begin{equation} \label{eq:xfiber}
 x  - \frac{1}{x} = s
\end{equation}
with
that of
\begin{equation} \label{eq:yfiber}
 y  + \frac{1}{y} = t - s.
\end{equation}
The set of solutions of (\ref{eq:xfiber}) degenerates
at $s = \pm 2 \sqrt{-1}$,
whereas that of (\ref{eq:yfiber}) degenerates
at $s = \pm 2 + t$.
% This shows that $W^{-1}(0)$ is a fourfold branched cover
% of the $s$-plane,
% whose branch points are $\pm 2$ and $\pm 2 \sqrt{-1}$.
The fiber $W^{-1}(0)$ is obtained
by gluing four copies $U_1, \ldots, U_4$
of the $s$-plane cut along four half lines $I_1, \ldots, I_4$
as in Figure \ref{fg:s-plane}.
The resulting surface is shown in Figure \ref{fg:P1P1_surface}.

Two of the four branch points of $\pi_t$ move
as in Figure \ref{fg:s-plane_motion}
as one varies $t$ from $0$ to
$p_1 = - 2 - 2 \sqrt{-1}$
along $c_1$.
This shows that
the image by $\pi_0$ of the corresponding vanishing cycle $C_1$
is the line segment from $2$ to $- 2 \sqrt{-1}$ up to homotopy.
Figure \ref{fg:s-plane_vc} shows four line segments on the $s$-plane,
which are images of the vanishing cycles of $W$
shown in Figure \ref{fg:P1P1_surface_cycles}.
One can equip these vanishing cycles with gradings
so that the Maslov indices of the intersection points in
$C_1 \cap C_2$,
$C_2 \cap C_3$, and
$C_3 \cap C_4$
are $1$.
Then the Maslov indices of
intersection points in $C_1 \cap C_4$ are $2$.
The dot in Figure \ref{fg:P1P1_surface_cycles} shows our choice
of branch points for the non-trivial spin structures.

\begin{figure}[htbp]
\begin{tabular}{cc}
\begin{minipage}{.5 \linewidth}
\centering
\input{P1P1_surface_graph.pst}
\caption{The graph on $W^{-1}(0)$}
\label{fg:P1P1_surface_graph}
\end{minipage}
&
\begin{minipage}{.5 \linewidth}
\centering
\input{P1P1_coamoeba.pst}
\caption{The coamoeba}
\label{fg:P1P1_coamoeba}
\end{minipage}
\ \\
\begin{minipage}{.5 \linewidth}
\centering
\input{P1P1_graph.pst}
\caption{The dimer model $G$}
\label{fg:P1P1_graph}
\end{minipage}
&
\begin{minipage}{.5 \linewidth}
\centering
\input{P1P1_matching.pst}
\caption{The perfect matching $D$}
\label{fg:P1P1_matching}
\end{minipage}
\end{tabular}
\end{figure}

There are four quadrangles in $W^{-1}(0)$
bounded by $\bigcup_{i=1}^4 C_i$,
which correspond to $A_\infty$-operations in $\Fuk W$.
By contracting these quadrangles,
one obtains the bicolored graph
in Figure \ref{fg:P1P1_surface_graph}.
Figure \ref{fg:P1P1_coamoeba} shows the coamoeba of $W^{-1}(0)$,
and Figure \ref{fg:P1P1_graph} shows the dimer model $G$.
Here, the black node on the top left in Figure \ref{fg:P1P1_graph}
corresponds to $U_1$ in Figure \ref{fg:P1P1_coamoeba}.

The perfect matching $D$
coming from the order in the distinguished basis of vanishing cycles
is shown in Figure \ref{fg:P1P1_matching}.
The full strong exceptional collection
of line bundles on $\bP^1 \times \bP^1$
associated with the pair $(G, D)$ is given by
\begin{equation} \label{eq:ec}
 (E_1, E_2, E_3, E_4) = (\scO, \scO(1,0), \scO(1,1), \scO(2,1)),
\end{equation}
where
$\scO(i,j) = \scO_{\bP^1}(i) \boxtimes \scO_{\bP^1}(j)$
denotes the exterior tensor product
of the $i$-th and $j$-th tensor powers
of the hyperplane bundle on $\bP^1$
for $i, j \in \bZ$.

\section{Non-isoradial dimer model and vanishing cycles}

In this section,
we discuss the role of the non-isoradial dimer model
in Figure \ref{fg:P1P1_graph2} below,
which corresponds to
the full strong exceptional collection
\begin{equation} \label{eq:ec2}
 (E_1, E_2, E_3, E_4) = (\scO, \scO(1,0), \scO(0,1), \scO(1,1))
\end{equation}
of line bundles on $\bP^1 \times \bP^1$.
This collection is obtained from the collection \eqref{eq:ec}
by the right mutation at the last term,
and the corresponding vanishing paths are shown
in Figure \ref{fg:P1P1_vp2}.
Figure \ref{fg:s-plane_vc2} shows the images of
vanishing cycles
by the projection to the $s$-plane.
The vanishing cycles on $W^{-1}(0)$ are shown 
in Figure \ref{fg:P1P1_surface_cycle2},
and by contracting eight triangles bounded by them,
one obtains the bicolored graph on $W^{-1}(0)$
shown in Figure \ref{fg:P1P1_surface_graph3}.

Unfortunately,
the argument projection of this graph
does not give a dimer model.
Figure \ref{fg:s-plane_vc3} shows the image of this graph
on the $s$-plane.
Here, two white dots are images of eight nodes,
and four black dots are branch points
where the fiber becomes two points
instead of four points.

\begin{figure}[htbp]
\begin{tabular}{cc}
\begin{minipage}{.5 \linewidth}
\centering
\input{P1P1_vp2.pst}
\caption{A distinguished set of vanishing paths}
\label{fg:P1P1_vp2}
\end{minipage}
&
\begin{minipage}{.5 \linewidth}
\centering
\input{s-plane_vc2.pst}
\caption{Images of vanishing cycles}
\label{fg:s-plane_vc2}
\end{minipage}
\\
\ \vspace{0mm}\\
\begin{minipage}{.5 \linewidth}
\centering
\input{P1P1_surface_cycle2.pst}
\caption{Vanishing cycles on $W^{-1}(0)$}
\label{fg:P1P1_surface_cycle2}
\end{minipage}
&
\begin{minipage}{.5 \linewidth}
\centering
\input{P1P1_surface_graph3.pst}
\caption{The graph on $W^{-1}(0)$}
\label{fg:P1P1_surface_graph3}
\end{minipage}
\end{tabular}
\end{figure}

\begin{figure}
\begin{tabular}{cc}
\begin{minipage}{.5 \linewidth}
\centering
\input{s-plane_vc3.pst}
\caption{The image of the graph on the $s$-plane}
\label{fg:s-plane_vc3}
\end{minipage}
&
\begin{minipage}{.5 \linewidth}
\centering
\input{s-plane_motion2.pst}
\caption{Behavior of the branch points}
\label{fg:s-plane_motion2}
\end{minipage}
\\
\begin{minipage}{.5 \linewidth}
\centering
\input{s-plane_vc4.pst}
\caption{The deformed graph}
\label{fg:s-plane_vc4}
\end{minipage}
&
\begin{minipage}{.5 \linewidth}
\centering
\input{s-plane_vc5.pst}
\caption{Labels on the graph on the $s$-plane}
\label{fg:s-plane_vc5}
\end{minipage}
\end{tabular}
\end{figure}

Consider the deformation
$$
 W_t(x, y) = x + \frac{2 t - 1}{x} + y + \frac{t + 1}{y}
$$
from
$$
 W_0(x, y) = x - \frac{1}{x} + y + \frac{1}{y}
$$
to
$$
 W_1(x, y) = x + \frac{1}{x} + y + \frac{2}{y}.
$$
Figure \ref{fg:s-plane_motion2} shows
the behavior of the branch points
under this deformation.
The graph in Figure \ref{fg:s-plane_vc3} is deformed
to the graph in Figure \ref{fg:s-plane_vc4}.
We put labels on this graph as in Figure \ref{fg:s-plane_vc5}.

Note that the map
$$
\begin{array}{ccc}
 \bCx & \to & \bC \\
 \vin & & \vin \\
  x & \mapsto & x + \frac{1}{x}
\end{array}
$$
from the $x$-plane to the $s$-plane is a branched double cover,
which maps both the upper half plane and the lower half plane
into the whole $s$-plane minus
$$
 I = I_+ \coprod I_- = \{ s \in \bR \mid |s| \ge 2 \},
$$
and the real axis into $I$.
The inverse image of the interval
$
 I_0 = \{ s \in \bR \mid |s| \le 2 \}
$
is the circle
$
 \{ x \in \bC \mid |x| = 1 \}.
$
Figures \ref{fg:s-plane_division} and
\ref{fg:x-plane_division} shows this behavior.
The same is true for the $s$-projection from the $y$-plane
$$
\begin{array}{ccc}
 \bCx & \to & \bC \\
 \vin & & \vin \\
  y & \mapsto & - y - \frac{2}{x}
\end{array}
$$
up to the minus sign and a rescaling.

\begin{figure}
\begin{tabular}{cc}
\begin{minipage}{.5 \linewidth}
\centering
\input{s-plane_division.pst}
\caption{Division of the $s$-plane}
\label{fg:s-plane_division}
\end{minipage}
&
\begin{minipage}{.5 \linewidth}
\centering
\input{x-plane_division.pst}
\caption{Division of the $x$-plane}
\label{fg:x-plane_division}
\end{minipage}
% \begin{minipage}{.5 \linewidth}
% \input{s-plane_vc6.pst}
% \caption{The graph on the $s$-plane}
% \label{fg:s-plane_vc6}
% \end{minipage}
% \begin{minipage}{.5 \linewidth}
% \input{s-plane_vc7.pst}
% \caption{The graph on the $s$-plane}
% \label{fg:s-plane_vc7}
% \end{minipage}
\end{tabular}
\end{figure}

\begin{figure}
\begin{tabular}{cc}
\begin{minipage}{.5 \linewidth}
\centering
\scalebox{.7}{\input{x-plane_graph.pst}}
\caption{The graph on the $x$-plane}
\label{fg:x-plane_graph}
\end{minipage}
&
\begin{minipage}{.5 \linewidth}
\centering
\scalebox{.7}{\input{y-plane_graph.pst}}
\caption{The graph on the $y$-plane}
\label{fg:y-plane_graph}
\end{minipage}
\end{tabular}
\end{figure}

\begin{figure}
\begin{tabular}{cc}
\begin{minipage}{.5 \linewidth}
\centering
\input{torus_graph.pst}
\caption{The graph on the torus}
\label{fg:torus_graph}
\end{minipage}
&
\begin{minipage}{.5 \linewidth}
\centering
\input{torus_graph_deform.pst}
\caption{A perturbed graph}
\label{fg:torus_graph_deform}
\end{minipage}
\\
\ \vspace{5mm}\\
\begin{minipage}{.5 \linewidth}
\centering
\input{torus_graph_deform2.pst}
\caption{A deformed graph}
\label{fg:torus_graph_deform2}
\end{minipage}
&
\begin{minipage}{.5 \linewidth}
\centering
\input{P1P1_graph2.pst}
\caption{The non-isoradial dimer model}
\label{fg:P1P1_graph2}
\end{minipage}
\end{tabular}
\end{figure}

Figures \ref{fg:x-plane_graph} and
\ref{fg:y-plane_graph} shows the pull-back of 
the graph on the $s$-plane in Figure \ref{fg:s-plane_vc5}
to the $x$-plane and the $y$-plane respectively.
This shows that the argument projection
of the graph on $W^{-1}(0)$
looks as in Figure \ref{fg:torus_graph}.
Although the restriction of the argument projection
to the graph is not injective,
one can perturb the graph slightly
to obtain the graph in Figure \ref{fg:torus_graph_deform}.
This graph is combinatorially identical to the graph
in Figure \ref{fg:torus_graph_deform2},
which gives the non-isoradial dimer model
in Figure \ref{fg:P1P1_graph2}.

{\bf Acknowledgment}:
K.~U. is supported by Grant-in-Aid for Young Scientists (No.18840029).

\bibliographystyle{plain}
\bibliography{bibs}

\def\cprime{$'$}
\begin{thebibliography}{1}

\bibitem{Feng-He-Kennaway-Vafa}
Bo~Feng, Yang-Hui He, Kristian~D. Kennaway, and Cumrun Vafa.
\newblock Dimer models from mirror symmetry and quivering amoebae.
\newblock {\em Adv. Theor. Math. Phys.}, 12(3):489--545, 2008.

\bibitem{Franco-Hanany-Martelli-Sparks-Vegh-Wecht_GTTGBT}
Sebasti{\'a}n Franco, Amihay Hanany, Dario Martelli, James Sparks, David Vegh,
  and Brian Wecht.
\newblock Gauge theories from toric geometry and brane tilings.
\newblock {\em J. High Energy Phys.}, (1):128, 40 pp. (electronic), 2006.

\bibitem{Franco-Hanany-Vegh-Wecht-Kennaway_BDQGT}
Sebasti{\'a}n Franco, Amihay Hanany, David Vegh, Brian Wecht, and Kristian~D.
  Kennaway.
\newblock Brane dimers and quiver gauge theories.
\newblock {\em J. High Energy Phys.}, (1):096, 48 pp. (electronic), 2006.

\bibitem{Hanany-Herzog-Vegh_BTEC}
Amihay Hanany, Christopher~P. Herzog, and David Vegh.
\newblock Brane tilings and exceptional collections.
\newblock {\em J. High Energy Phys.}, (7):001, 44 pp. (electronic), 2006.

\bibitem{Hanany-Kennaway_DMTD}
Amihay Hanany and Kristian~D. Kennaway.
\newblock Dimer models and toric diagrams.
\newblock hep-th/0503149, 2005.

\bibitem{Hanany-Vegh}
Amihay Hanany and David Vegh.
\newblock Quivers, tilings, branes and rhombi.
\newblock {\em J. High Energy Phys.}, (10):029, 35, 2007.

\bibitem{Seidel_VC2}
Paul Seidel.
\newblock More about vanishing cycles and mutation.
\newblock In {\em Symplectic geometry and mirror symmetry (Seoul, 2000)}, pages
  429--465. World Sci. Publishing, River Edge, NJ, 2001.

\bibitem{Ueda-Yamazaki_NBTMQ}
Kazushi Ueda and Masahito Yamazaki.
\newblock A note on dimer models and {M}c{K}ay quivers.
\newblock math.AG/0605780.

\bibitem{Ueda-Yamazaki_toricdP}
Kazushi Ueda and Masahito Yamazaki.
\newblock Homological mirror symmetry for toric orbifolds of toric del {P}ezzo
  surfaces.
\newblock math.AG/0703267.

\end{thebibliography}

\noindent
Kazushi Ueda

Department of Mathematics,
Graduate School of Science,
Osaka University,
Machikaneyama 1-1,
Toyonaka,
Osaka,
560-0043,
Japan.

{\em e-mail address}\ : \  kazushi@math.sci.osaka-u.ac.jp

\ \\

\noindent
Masahito Yamazaki

Department of Physics,
Graduate School of Science,
University of Tokyo,
Hongo 7-3-1,
Bunkyo-ku,
Tokyo,
113-0033,
Japan

{\em e-mail address}\ : \  yamazaki@hep-th.phys.s.u-tokyo.ac.jp

\end{document}